
\input amstex.tex
\documentstyle{amsppt}
\magnification1100

\hsize=12.5cm
\vsize=17cm
\hoffset=1cm
\voffset=2cm

\footline={\hss{\vbox to 2cm{\vfil\hbox{\rm\folio}}}\hss}
\nopagenumbers
\def\DJ{\leavevmode\setbox0=\hbox{D}\kern0pt\rlap
{\kern.04em\raise.188\ht0\hbox{-}}D}

\baselineskip=13pt
\def\hf{{\textstyle{1\over2}}}
\def\a{\alpha}
\def\d{{\,\roman d}}
\def\e{\varepsilon}

\def\G{\Gamma}

\def\s{\sigma}

\def\={\;=\;}

\def\zt{\zeta(\hf+it)}

\def\D{\Delta}

\def\R{\Re{\roman e}\,} 
\def\z{\zeta}

\def\hf{{\textstyle{1\over2}}}

\font\tenmsb=msbm10
\font\sevenmsb=msbm7
\font\fivemsb=msbm5
\newfam\msbfam
\textfont\msbfam=\tenmsb
\scriptfont\msbfam=\sevenmsb
\scriptscriptfont\msbfam=\fivemsb
\def\Bbb#1{{\fam\msbfam #1}}

\def \NN {\Bbb N}

\def \RR {\Bbb R}

\font\ff=cmr8

\baselineskip=13pt

\font\teneufm=eufm10
\font\seveneufm=eufm7
\font\fiveeufm=eufm5
\newfam\eufmfam
\textfont\eufmfam=\teneufm
\scriptfont\eufmfam=\seveneufm
\scriptscriptfont\eufmfam=\fiveeufm
\def\mathfrak#1{{\fam\eufmfam\relax#1}}

\font\tenmsb=msbm10
\font\sevenmsb=msbm7
\font\fivemsb=msbm5
\newfam\msbfam
     \textfont\msbfam=\tenmsb
      \scriptfont\msbfam=\sevenmsb
      \scriptscriptfont\msbfam=\fivemsb

  \def\rightheadline{{\hfil{\ff
 Some remarks on the moments of $|\zt|$ in short intervals}\hfil\tenrm\folio}}

  \def\leftheadline{{\tenrm\folio\hfil{\ff
   Aleksandar Ivi\'c }\hfil}}
  \def\emptyheadline{\hfil}
  \headline{\ifnum\pageno=1 \emptyheadline\else
  \ifodd\pageno \rightheadline \else \leftheadline\fi\fi}

\topmatter
\title
Some remarks on the moments of $|\zt|$ in short intervals
\endtitle
\author   Aleksandar Ivi\'c  \endauthor
\address
Aleksandar Ivi\'c, Katedra Matematike RGF-a,
Universitet u Beogradu, \DJ u\v sina 7, 11000 Beograd,
Serbia.
\endaddress
\email {\tt
ivic\@rgf.bg.ac.yu,  aivic\@matf.bg.ac.yu} \endemail
\dedicatory
\enddedicatory
\abstract
{Some new results on power moments of the integral
$$
J_k(t,G) = {1\over\sqrt{\pi}G}
\int_{-\infty}^\infty |\z(\hf + it + iu)|^{2k}{\roman e}^{-(u/G)^2}\d u
\qquad(t \asymp T, T^\e \le G \ll T,\,k\in\NN)
$$
are obtained when $k=1$.
These results can be used to derive
bounds for moments of $|\zt|$.}
\endabstract
\endtopmatter
\document

\head
1. Introduction
\endhead
Power moments represent one of the most important parts of the
theory of the Riemann zeta-function $\z(s) = \sum_{n=1}^\infty n^{-s}$
$\;(\s = \R s > 1)$. Of particular significance are the moments on the
``critical line" $\s = \hf$, and a vast literature exists on this
subject (see e.g., the monographs [2], [3], [12]). Let us define
$$
I_k(T) = \int_0^T|\zt|^{2k}\d t,\leqno(1.1)
$$
where $k\in\RR$ is a fixed, positive number.
The aim of this paper is to investigate upper bounds for $I_k(T)$
when $k\in\NN$, which we henceforth assume. The problem can be reduced
to bounds of $|\zt|$ over short intervals,  but it
is more expedient to work with the smoothed integral
$$
J_k(T,G) := {1\over \sqrt{\pi}G}\int_{-\infty}^\infty
|\z(\hf + iT + iu)|^{2k}{\roman e}^{-(u/G)^2}\d u\quad(1 \ll G \ll T).
\leqno(1.2)
$$
Namely we obviously have
$$
I_k(T+G) - I_k(T-G) = \int_{-G}^G|\z(\hf + iT + iu)|^{2k}\d u
\le \sqrt{\pi}{\roman e}G\,J_k(T,G),\leqno(1.3)
$$
and it is technically more convenient to work with $J_k(T,G)$ than with
the differenced integral
$I_k(T+G) - I_k(T-G)$. Of course, instead of the Gaussian exponential
weight $\exp(-(u/G)^2)$, one could introduce in (1.2) other smooth weights
with a similar effect. The Gaussian weight has the advantage that, by the
use of the classical integral
$$
\int_{-\infty}^\infty \exp(Ax-Bx^2)\d x = \sqrt{\pi\over B}\,\exp\left(
{A^2\over4B}\right)\qquad(\R B > 0),\leqno(1.4)
$$
one can often explicitly evaluate the relevant exponential integrals
that appear in the course of the proof.

\medskip
One expects that $J_1(t,G)$, at least for certain ranges of $G= G(T)$, behaves
in $\,[T,\,2T]$ like $O(t^\e)$ on the average. This would be a trivial consequence,
for $1\ll G\ll T$, of the truth of the famous Lindel\"of hypothesis that
$\zt \ll_\e |t|^\e$.
In [5] we proved the following result on moments
of $J_1(t,G)$ which supports this claim. Our bounds were given by

\bigskip
THEOREM A. {\it We have}
$$
\int_T^{2T}J_1^m(t,G)\d t \ll_\e T^{1+\e}\leqno(1.5)
$$
{\it for $T^\e \le G \le T$ if $\,m = 1,2$; for $T^{1/7+\e} \le G \le T$ if
$\,m=3$, and for $T^{1/5+\e} \le G \le T$ if $m=4$.}

\bigskip
 Here and later $\e >0$ denotes
constants which may be arbitrarily small, but are not necessarily the
same ones at each occurrence, while $a \ll_\e b$ means that the
$\ll$--constant depends only on $\e$. It is the lower bound for $G$ in the above
results that matters, because for $T^{1/3} \le G = G(T) \ll T$ the bound in (1.5)
trivially holds, since (see [2, Chapter 7])
$$
\int_{T-G}^{T+G}|\zt|^2\d t \ll G\log T\qquad(T^{1/3} \le G \ll T).
$$

In the case when $m = 4$  we can improve on the range of $G$
furnished by Theorem A,
and when  $m = 5$ and $m=6$ we can obtain  new results. This is given by

\medskip
THEOREM 1. {\it We have } (1.5) {\it with $T^{7/36} \le G = G(T) \le T$ when $m = 4$,
with $T^{1/5} \le G = G(T) \le T$ when $m=5$, and with $T^{2/9} \le G = G(T) \le T$
when $m=6$.}

\medskip

While the proof of Theorem A in [5] rested on the explicit formula for
$J_1(T,G)$ (cf. Lemma 2) and direct evaluation, in obtaining the
improvement contained in Theorem 1 we shall deal with
the moments of the error term function
$$
E^*(t) \;:=\; E(t) - 2\pi\D^*\bigl({t\over2\pi}\bigr),\leqno(1.6)
$$
where as usual ($\gamma = -\G'(1) = 0.5772157\dots\,$ is the Euler constant)
$$
E(T) \;=\;\int_0^T|\zt|^2\d t - T\left(\log\bigl({T\over2\pi}\bigr) + 2\gamma - 1
\right)\leqno(1.7)
$$
is the error term in the mean square formula for $|\zt|$ (see [2] and [3] for
a comprehensive account). Furthermore
$$
\D^*(x) := -\D(x)  + 2\D(2x) - \hf\D(4x)
= \hf\sum_{n\le4x}(-1)^nd(n) - x(\log x + 2\gamma - 1),
\leqno(1.8)
$$
with
$$
\D(x) \;=\; \sum_{n\le x}d(n) - x(\log x + 2\gamma - 1)
\leqno(1.9)
$$
the error term in the classical Dirichlet divisor problem, and $d(n)$ the
number of divisors of $n$. The function
$E^*(t)$ is smaller on the average than $E(t)$, and for this reason
it seems better suited to use it in our context than $E(t)$.

\medskip
Theorem 1 follows from the bounds for moments of $E^*(t)$ and

\medskip
THEOREM 2. {\it For $\,T^\e \ll G = G(T) \le T$ and fixed $m\ge 1$ we have}
$$
\int_T^{2T}J_1^m(t,G)\d t \ll G^{-1-m}\int_{-G\log T}^{G\log T}
\left(\int_T^{2T}|E^*(t+x)|^m\d t\right)\d x + T\log^{2m}T.\leqno(1.10)
$$

\medskip

\medskip
The plan of the paper is as follows. In Section 2 we shall present the
results on the moments of $E^*(t)$. Then, in Section 3, we shall prove
both Theorem 1 and Theorem 2. Some concluding remarks will be given
in Section 4.

\head
2. The moments of $E^*(t)$
\endhead
M. Jutila [7], [8] investigated both the
local and global behaviour of the function $E^*(t)$,
and in particular in [8] he proved that
$$
\int_0^T(E^*(t))^2\d t \;\ll\; T^{4/3}\log^3T.\leqno(2.1)
$$
This bound is remarkable, because (see  [3, Theorem 2.4])
$$
\int_0^T E^2(t)\d t = cT^{3/2} + O(T\log^5T),\quad
c = {2\over3}(2\pi)^{-1/2}\,{\z^4(3/2)\over\z(3)} = 10.3047\ldots\,,
$$
which shows that, in the mean square sense, the function $E^*(t)$
is much smaller than $E(t)$ or, in other words, the functions
$E(t)$ and $2\pi\D^*(t/(2\pi))$ are ``close" to one another.

\medskip
In the first part of the author's work [4] the bound in (2.1) was complemented
with the new bound
$$
\int_0^T (E^*(t))^4\d t \;\ll_\e\; T^{16/9+\e};\leqno(2.2)
$$
neither (2.1) or (2.2) seem to imply each other.
In the second part of the same work (op. cit.) it was proved that
$$
\int_0^T |E^*(t)|^5\d t \;\ll_\e\; T^{2+\e},\leqno(2.3)
$$
and some further results on higher moments of $|E^*(t)|$ were obtained as well.
In [6] the author sharpened (2.1) by proving the asymptotic formula
$$
\int_0^T (E^*(t))^2\d t \;=\; T^{4/3}P_3(\log T) + O_\e(T^{7/6+\e}),\leqno(2.4)
$$
where $P_3(y)$ is a polynomial of degree three in $y$ with
positive leading coefficient, and all the coefficients may be evaluated
explicitly.

\medskip
In the third part of [5] the integral of $E^*(t)$ was investigated.
If  we define the error-term function $R(T)$ by the relation
$$
\int_0^T E^*(t)\d t = {3\pi\over 4}T + R(T),\leqno(2.5)
$$
then
we have (see  [1], [3] for the first formula and [2] for the second one)
$$
\int_0^T E(t)\d t = \pi T + G(T),\quad \int_0^T \D(t)\d t = {T\over 4} + H(T),
$$
where both $G(T), H(T)$ are $O(T^{3/4})$ and also $\Omega_\pm(T^{3/4})$
(for $g(x) > 0\;(x>x_0)\,$ $f(x) = \Omega(g(x))$ means that $f(x) = o(g(x))$
does not hold as $x\to\infty$, $f(x) = \Omega_\pm(g(x)))$ means that there
are unbounded sequences $\{x_n\},\,\{y_n\},\,$ and constants $A,B>0$ such
that $f(x_n) > Ag(x_n)$ and $f(y_n) < -Bg(y_n)$). Since
$$
\int_0^T \D(at)\d t = {1\over a}\int_0^{aT}\D(x)\d x
\qquad(a>0,\;T>0)
$$
holds, it is obvious that ${3\pi\over 4}$ is the
``correct" constant in (2.5), and that trivially one has the bound
$R(T) = O(T^{3/4})$. We  proved in [4, Part III] the following results:
\bigskip
THEOREM B. {\it We have}
$$
R(T) \;=\; O_\e(T^{593/912+\e}), \quad {593\over912} = 0.6502129\ldots\,.
\leqno(2.6)
$$

\bigskip
THEOREM C. {\it We have
$$
\int_0^TR^2(t)\d t \;=\; T^2p_3(\log T) + O_\e(T^{11/6+\e}),\leqno(2.7)
$$
where $p_3(y)$ is a cubic polynomial in $y$ with positive leading
coefficient, whose all coefficients may be explicitly evaluated.}

\bigskip
THEOREM D. {\it We have}
$$
\int_0^TR^4(t)\d t \;\ll_\e\; T^{3+\e}.
\leqno(2.8)
$$

\bigskip
These results shows that $E^*(t)$ and its integral are smaller
on the average than $E(t)$ and its integral, respectively (see [2]
and [3]). Thus the effect of Theorem 2 is that the moments of $J_1$
are bounded by moments of $E^*$, and not by the moments of $E$ itself.

\medskip

Note that (2.4) and (2.7) imply, respectively,
$$
E^*(T) = \Omega(T^{1/6}(\log T)^{3/2}),\quad
R(T) = \Omega(T^{1/2}(\log T)^{3/2}),
$$
while we have
$$
E(T) =\Omega(T^{1/4}L(T)),\quad \D^*(x) = \Omega(x^{1/4}L(x)),\leqno(2.9)
$$
where ($\log_kx = \log(\log_{k-1}x))$ for $k\ge 2$, $\log_1x \equiv\log x$)
$$
L(y) := (\log y)^{1/4}(\log_2y)^{{3\over4}(2^{4/3}-1)}(\log_3 y)^{-5/8}
\qquad(y > {\roman e}^{\roman e}).
\leqno(2.10)
$$
These are the strongest known $\Omega$--results for $E(T)$ and $\D^*(x)$,
and follow by the method of K. Soundararajan [11], who obtained the analogue
of (2.9)--(2.10) for $\D(x)$. Lau--Tsang [9] proved the first $\Omega$--result
in (2.9), and the second one follows by their Theorem 1.2 and Soundararajan's
result for $\D(x)$.
\head
3. The  proof of theorem 1 and Theorem 2
\endhead
We suppose that $1\ll G = G(T) \le T^{1/3}, T\le t\le 2T$. From (1.2)
we have, on integrating by parts,
$$
\eqalign{
J_1(t,G) &= {1\over\sqrt{\pi}G}\int_{-\infty}^\infty \Bigl(\log(t+x)+2\gamma
- \log(2\pi) + E'(t+x)\Bigr){\roman e}^{-(x/G)^2}\d x\cr&
= O(\log T) + {2\over\sqrt{\pi}G^3}\int_{-\infty}^\infty
xE(t+x){\roman e}^{-(x/G)^2}\d x,\cr}\leqno(3.1)
$$
and we may truncate the last integral at $x = \pm GL$, with a very
small error, where for brevity we put $L :=\log T$. But from (1.6)  we have
$$
\eqalign{&
\int_{-GL}^{GL} xE(t+x){\roman e}^{-(x/G)^2}\d x
 = \int_{-GL}^{GL}  xE^*(t+x){\roman e}^{-(x/G)^2}\d x\cr&
+2\pi\int_{-GL}^{GL}  x\D^*\left({t+x\over2\pi}\right){\roman e}^{-(x/G)^2}\d x.
\cr}\leqno(3.2)
$$
From (1.8) it follows, on integrating by parts, that
$$
\eqalign{&
\int_{-GL}^{GL}  x\D^*(t+x){\roman e}^{-(x/G)^2}\d x\cr&
= 2\int_{-GL}^{GL}  x\sum_{n\le2(t+x)/\pi}(-1)^nd(n)\cdot{\roman e}^{-(x/G)^2}\d x
+ O(G^3\log T)\cr&
= \Bigl({\pi\over2}\Bigr)^2
\int_{-GL}^{GL}  y\sum_{n\le{2t\over\pi}+y}(-1)^nd(n)\cdot{\roman e}^{-(\pi y/(2G))^2}
\d y + O(G^3\log T)\cr&
= \hf G^2\int_{-GL}^{GL}{\roman e}^{-(\pi y/(2G))^2}\d\left(
\sum_{n\le{2t\over\pi}+y}(-1)^nd(n)\right) + O(G^3\log T)\cr&
\ll G^2\sum_{{2t\over\pi}-GL\le n \le {2t\over\pi}+GL}d(n) + G^3\log T
\ll G^3L\log T,\cr}
$$
on using a result of P. Shiu [10] on multiplicative functions in short intervals.
Therefore from (3.1)--(3.2) we obtain ($T^\e \le G = G(T) \le T^{1/3},\,T\le t \le 2T$)
$$
J_1(t,G) = {2\over\sqrt{\pi}G^3}\int_{-G\log T}^{G\log T}
xE^*(t+x){\roman e}^{-(x/G)^2}\d x + O(\log^2 T).\leqno(3.3)
$$
By H\"older's inequality it follows that, for $m > 1$ fixed
(not necessarily an integer; for $m=1$ the assertion is easy)
$$
\eqalign{&
\int_T^{2T}J_1^m(t,G)\d t \ll {1\over G^{3m}}\int_T^{2T}\left(
\int_{-GL}^{GL}|xE^*(t+x)|{\roman e}^{-(x/G)^2}\d x\right)^m\d t + T\log^{2m}T\cr&
\ll {1\over G^{3m}}\int_T^{2T}\int_{-GL}^{GL}|E^*(t+x)|^m\d x {\left(\int_{0}^{GL}
x^{m\over m-1}{\roman e}^{-(x/G)^2}\d x\right)}^{m-1}\d t + T\log^{2m}T\cr&
\ll G^{-1-m}\int_{-GL}^{GL}\left(\int_T^{2T}|E^*(t+x)|^m\d t\right)\d x
+ T\log^{2m}T,
\cr}
$$
as asserted by (1.10) of Theorem 2.

\medskip
To prove Theorem 1, we use first (2.2) and Theorem 2 with $m=4$. We obtain
$$
\int_T^{2T}J_1^4(t,G)\d t \ll_\e G^{-5}GT^{16/9+\e} + T\log^{10}T \ll_\e T^{1+\e}
$$
for $G \ge T^{7/36}$. Similarly from (2.3) and Theorem 2 with $m=5$ we have
$$
\int_T^{2T}J_1^5(t,G)\d t \ll_\e G^{-6}GT^{2+\e} + T\log^5T \ll_\e T^{1+\e}
$$
for $G \ge T^{1/5}$.

\smallskip
It remains to deal with the case $m=6$, when the assertion of Theorem 1 follows
analogously from Theorem 2 and the bound
$$
\int_0^T|E^*(t)|^6\d t \;\ll_\e\; T^{7/3+\e}.\leqno(3.4)
$$
To obtain (3.4) note that
$$
\eqalign{&
\int_T^{2T}|E^*(t)|^6\d t \,=\;  \int_{|E^*|\le T^{1/3}}
+ \;\;\int_{|E^*|> T^{1/3}}\cr&
\,\le T^{1/3}\int_T^{2T}|E^*(t)|^5\d t
+ T^{-2/3}\int_T^{2T}|E^*(t)|^8\d t\cr&
\,\ll_\e \;T^{7/3+\e}.\cr}
$$
Here we used (2.3), (1.6) and (see [3])
$$
\int_T^{2T}|E^*(t)|^8\d t \ll \int_T^{2T}|E(t)|^8\d t
+ \int_1^{20T}|\D(t)|^8\d t \ll_\e T^{3+\e}.
$$
This completes the proof of Theorem 1, with the remark that
bounds for higher moments of $|E^*(t)|$ could be also
derived, but their sharpness would decrease as $m$ increases.

\head
4. Concluding remarks
\endhead

In [5] the author proved the following result, which connects
the moments of $|\zt|$ to the moments of $J_k(t,G)$. This is

\medskip
THEOREM E. {\it Suppose that}
$$
\int_T^{2T}J_k^m(t,G)\d t \;\ll_\e\; T^{1+\e}
$$
{\it holds for some fixed $k,m  \in\NN$ and $G = G(T) \ge T^{\a_{k,m}+\e},\,
 0 \le \a_{k,m} < 1$. Then}
$$
\int_0^T|\zt|^{2km}\d t \;\ll_\e\; T^{1+(m-1)\a_{k,m}+\e}.\leqno(4.1)
$$
If we use (4.1) with $k=1$ together with (1.10) of Theorem 2, then
we obtain that bounds for moments of $|\zt|$ can be found directly from
the bounds for moments
of $E^*(t)$. From Theorem 1 it follows that we can take $\a_{1,4} = 7/36,
\,\a_{1,5} = 1/5$, $\a_{1,6} = 2/9$. Therefore we obtain from (4.1)
$$\eqalign{&
\int_0^T|\zt|^{8}\d t \ll_\e T^{19/12+\e},
\;
\int_0^T|\zt|^{10}\d t \ll_\e T^{9/5+\e},\cr&
\int_0^T|\zt|^{12}\d t \ll_\e T^{19/9+\e}.\cr}\leqno(4.2)
$$
The exponents in (4.2) are somewhat poorer than the best known exponents (see
[3, Chapter 8]) 3/2, 7/4 and 2, respectively. However, it is clear that moments
of $J_k(t,G)$ are important for the estimation of moments of $|\zt|$, one of
central topics in zeta-function theory.

\medskip
One can obtain even a more direct connection between the moments of $|\zt|$ and
$|E^*(t)|$, as was shown in [4, Part II]. This is

\bigskip
THEOREM F.
{\it Let $k \ge 1$ be a fixed real, and let $c(k)$ be such a
constant for which
$$
\int_0^T|E^*(t)|^{k}\d t \;\ll_\e\; T^{c(k)+\e}.\leqno(4.3)
$$
Then we have
$$
\int_0^T|\zt|^{2k+2}\d t \;\ll_\e\; T^{c(k)+\e}.\leqno(4.4)
$$}

From (2.4) it follows that one can take $c(5) = 2$ in (4.3), so that (4.4)
gives the estimate
$$
\int_0^T|\zt|^{12}\d t \;\ll_\e\; T^{2+\e},\leqno(4.5)
$$
which is a result of D.R.
Heath-Brown [2], and it is (up to `$\e$') the strongest known bound of its kind.
With $c(6) = 7/3$, which follows from (3.4), we obtain from (4.4)
$$
\int_0^T|\zt|^{14}\d t \;\ll_\e\; T^{7/3+\e},
$$
but this follows trivially from (4.5) and the classical bound $\zt \ll |t|^{1/6}$.
\vfill
\eject
\topglue1cm
\bigskip
\Refs
\bigskip

\item{[1]} J.L. Hafner and A. Ivi\'c, On the mean square of the
Riemann zeta-function on the critical line,  J. Number Theory,
{\bf 32}(1989), 151--191.

\item{[2]} A. Ivi\'c, The Riemann zeta-function, John Wiley \&
Sons, New York, 1985.

\item{[3]} A. Ivi\'c, The mean values of the Riemann zeta-function,
LNs {\bf 82}, Tata Inst. of Fundamental Research, Bombay (distr. by
Springer Verlag, Berlin etc.), 1991.

\item{[4]} A. Ivi\'c, On the Riemann zeta-function and the divisor problem,
Central European J. Math. {\bf(2)(4)} (2004), 1-15, and II, ibid.
{\bf(3)(2)} (2005), 203-214,
and III, subm. to Ann. Univ. Budapest. Sectio Computatorica.

\item{[5]} A. Ivi\'c, On moments of $|\zt|$ in short intervals,
Proc. Conference in Honour of K. Ramachandra, Chennai, 2003,
Ramanujan Math. Soc., 2006, in press.

\item{[6]} A. Ivi\'c, On the mean square of the zeta-function and
the divisor problem, Annales  Acad. Scien. Fennicae Mathematica, in press.

\item{[7]} M. Jutila, Riemann's zeta-function and the divisor problem,
Arkiv Mat. {\bf21}(1983), 75-96 and II, ibid. {\bf31}(1993), 61-70.

\item{[8]} M. Jutila, On a formula of Atkinson, in ``Coll. Math. Sci.
J\'anos Bolyai 34, Topics in classical Number Theory, Budapest 1981",
North-Holland, Amsterdam, 1984, pp. 807-823.

\item{[9]} Y.-K. Lau and K.-M. Tsang, Omega result for the mean square of the
Riemann zeta-function, Manuscripta Math.  {\bf115}(2005), 373-381.

\item{[10]} P. Shiu, A Brun--Titchmarsh theorem for multiplicative
functions, J. reine angew. Math. {\bf31}(1980), 161-170.

\item{[11]} K. Soundararajan, Omega results for the divisor and circle problems,
Int. Math. Research Notices  {\bf36}(2003), 1987-1998.

\item{[12]} E.C. Titchmarsh, The theory of the Riemann zeta-function
(2nd ed.),  University Press, Oxford, 1986.

\endRefs

\enddocument

\bye